\documentclass[
a4paper,
11pt,
twoside,
english,
]{article}

\usepackage[utf8]{inputenc}
\usepackage[T1]{fontenc}
\usepackage[sc,noBBpl]{mathpazo}
\usepackage{babel}
\usepackage{geometry}
\usepackage{amsmath,amsfonts,amssymb,amsthm}
\usepackage{graphicx}
\usepackage{mathtools}
\usepackage[dvipsnames,table]{xcolor}
\usepackage{pdflscape}
\usepackage{etoolbox}
\usepackage{fancyhdr}
\usepackage{hyphenat}
\usepackage{enumitem}
\usepackage{siunitx}
\usepackage{booktabs}

\theoremstyle{plain}

\newcommand{\dd}{\mathop{}\!\mathrm{d}}

\usepackage[
  mincrossrefs=99,
  style=numeric-comp,
  bibstyle=siam,
  sorting=nyt,
  backref=true,
  backend=biber,
  url=true,
  doi=true,
  maxbibnames=99,
  maxnames=99,
  giveninits=true,
  useprefix=true,
  date=year,
  autolang=other
]{biblatex}
\DeclareCiteCommand{\citeauthornsc}
  {%
   \boolfalse{citetracker}%
   \boolfalse{pagetracker}%
   \usebibmacro{prenote}}
  {\ifciteindex
     {\indexnames{labelname}}
     {}%
   \printnames{labelname}}
  {\multicitedelim}
  {\usebibmacro{postnote}}
\DeclareSourcemap{
  \maps[datatype=bibtex]{
    \map{
      \step[fieldsource=doi,final]
      \step[fieldset=url,null]
    }
  }
}
\DeclareNumChars*{A}

\usepackage{csquotes}
\usepackage{fnpct}

\definecolor{webgreen}{rgb}{0,.5,0}
\definecolor{webbrown}{rgb}{.6,0,0}
\definecolor{myblue}{rgb}{0,0.25,0.5}
\usepackage[
  colorlinks=true,
  breaklinks=true,
  bookmarksnumbered,
  pdfhighlight=/O,
  urlcolor=webbrown,
  linkcolor=myblue,
  citecolor=webgreen,
  hyperfootnotes=true,
]{hyperref}

\newcommand{\email}[1]{\href{mailto:#1}{\texttt{#1}}}

\numberwithin{equation}{section}
\numberwithin{table}{section}
\numberwithin{figure}{section}

\makeatletter
\newcommand{\pagerefstar}{\@pagerefstar}
\makeatother

\usepackage{caption}
\newcommand{\mydate}{March 8, 2022}

\fancypagestyle{mypagestyle}{%
  \fancyhf{}%
  \fancyhead[LO,RE]{\scriptsize Pseudospectral method for population models with two structures}%
  \fancyhead[RO,LE]{\scriptsize Andò, De Reggi, Liessi, Scarabel}%
  \fancyfoot[LE,RO]{\scriptsize \thepage\,/\,\pagerefstar{LastPage}}%
  \fancyfoot[LO,RE]{\scriptsize \mydate{}}%
}
\usepackage{lastpage}

\allowdisplaybreaks

\begin{document}
\thispagestyle{empty}

\begin{center}
\LARGE
A pseudospectral method for investigating the stability of linear population models with two physiological structures

\bigskip
\large
Alessia Andò\footnotemark[1]\footnotemark[4]\footnotemark[6],
Simone De Reggi\footnotemark[2]\footnotemark[4]\footnotemark[7],
\\
Davide Liessi\footnotemark[2]\footnotemark[4]\footnotemark[8],
Francesca Scarabel\footnotemark[3]\footnotemark[4]\footnotemark[5]\footnotemark[9]

\footnotetext[1]{Area of Mathematics, Gran Sasso Science Institute, Viale F. Crispi 7, 67100 L'Aquila, Italy}%
\footnotetext[2]{Department of Mathematics, Computer Science and Physics, University of Udine, Via delle Scienze 206, 33100 Udine, Italy}%
\footnotetext[3]{Department of Mathematics, The University of Manchester, Oxford Rd, M13 9PL, Manchester, United Kingdom}%
\footnotetext[4]{CDLab -- Computational Dynamics Laboratory, University of Udine, Italy}%
\footnotetext[5]{Joint UNIversities Pandemic and Epidemiological Research, United Kingdom}%
\footnotetext[6]{\email{alessia.ando@gssi.it}}%
\footnotetext[7]{\email{dereggi.simone@spes.uniud.it}}%
\footnotetext[8]{\email{davide.liessi@uniud.it}}%
\footnotetext[9]{\email{francesca.scarabel@manchester.ac.uk}}%
\setcounter{footnote}{9}%

\medskip
\mydate
\end{center}

\begin{abstract}
The asymptotic stability of the null equilibrium of a linear population model with two physiological structures formulated as a first-order hyperbolic PDE is determined by the spectrum of its infinitesimal generator.
We propose an equivalent reformulation of the problem in the space of absolutely continuous functions in the sense of Carathéodory, so that the domain of the corresponding infinitesimal generator is defined by trivial boundary conditions.
Via bivariate collocation, we discretize the reformulated operator as a finite-dimensional matrix, which can be used to approximate the spectrum of the original infinitesimal generator.
Finally, we provide test examples illustrating the converging behavior of the approximated eigenvalues and eigenfunctions, and its dependence on the regularity of the model coefficients.


\end{abstract}


\section{Introduction}

Mathematical models are often used to describe the evolution of populations in biology and epidemiology. An important class of models that has attracted increased attention is that of structured population models, in which individuals are characterized by one or more variables that describe the i-state (i.e., the individual state) and determine the individual processes, including for instance birth, growth and death.
Examples of physiological structures are age, size, spatial position or time since infection.
We here focus on a class of structured population models where the structuring variables are continuous and the models are formulated as first-order hyperbolic partial differential equations (PDEs) (see, e.g., \cite{BekkalBrikciClairambaultRibbaPerthame2008,DysonVillellaBressanWebb2000,Howard2003,MagalRuan2008,MetzDiekmann1986} and references therein).
In particular, we consider the case where individuals are characterized by two different traits \cite{KangHuoRuan2021,Webb1985c}.

In applications to population dynamics, the interest is often focused on the long-term properties of the systems, for instance the existence of equilibrium states and their stability.
For linear models with two structures, it has been proved that the stability of the zero equilibrium is determined by the spectrum of the infinitesimal generator (IG) of the semigroup of solution operators \cite{KangHuoRuan2021,Webb1985c}.
Since the IG is an operator acting on an infinite-dimensional space of functions, numerical techniques are required to obtain finite-dimensional approximations of the operator and, in turn, of its spectrum.

For the analysis of local stability of equilibria, pseudospectral methods have been widely used both for delay equations \cite{BredaDiekmannGyllenbergScarabelVermiglio2016,BredaDiekmannMasetVermiglio2013,BredaGettoSanchezSanzVermiglio2015,BredaMasetVermiglio2005} and for PDE population models with one structuring variable \cite{BredaCusulinIannelliMasetVermiglio2007,BredaIannelliMasetVermiglio2008,ScarabelBredaDiekmannGyllenbergVermiglio2021}. The main advantage of pseudospectral methods is their typical spectral accuracy, by which the order of convergence of the approximation error increases with the regularity of the approximated function. In particular, the convergence is exponential for analytic functions \cite{Trefethen2013}. In the case of delay equations, this implies that the spectrum of the IG is approximated with exponential order of convergence, as the corresponding eigenfunctions are exponentials (see, e.g., \cite[Proposition 3.4]{BredaMasetVermiglio2015}). In the case of PDEs, the eigenfunctions are still exponential in time, but the order of regularity with respect to the physiological variables depends on the regularity of the model parameters, which therefore affects the order of convergence of the approximation \cite[Theorem 3.2]{Webb1985c}.
A similar behavior has been shown for the approximation of $R_0$ for structured epidemic models \cite{BredaDeReggiScarabelVermiglioWu2021,BredaFlorianRipollVermiglio2021,BredaKuniyaRipollVermiglio2020}, and in the approximation of the solution operators \cite{BredaLiessi2018,BredaLiessi2020,BredaMasetVermiglio2012}.

For structured models with one single structuring variable, pseudospectral methods have already been proposed to study the stability of the zero equilibrium in \cite{BredaCusulinIannelliMasetVermiglio2007}.
In that paper the IG is approximated by combining pseudospectral differentiation with the inversion of a (linear) algebraic condition characterizing the domain of the operator.
However, with more structuring variables implementing this technique becomes substantially more involved.

A different approach, which has been successfully employed in the context of nonlinear PDEs with one physiological variable \cite{ScarabelBredaDiekmannGyllenbergVermiglio2021} and of renewal equations \cite{ScarabelDiekmannVermiglio2021,ScarabelVermiglio}, consists in first reformulating the problem at hand via conjugation with an integral operator, and then approximating the resulting transformed operator via pseudospectral (also known as spectral collocation) techniques.
The advantages of this approach are mainly twofold: on the one hand, the transformed operator acts on a space of absolutely continuous functions (rather than the original space $L^1$), hence point evaluation, as well as polynomial interpolation and collocation, are well defined; on the other hand, the domain of the transformed operator is characterized by a trivial condition (specifically, a zero boundary condition), which substantially simplifies the numerical implementation.
From a modeling point of view, the integrated state has a clear interpretation as it represents the number of individuals whose i-state is less than a given value.

Goal of this work is to introduce a numerical method for the stability analysis of linear PDE population models with two structuring variables based on this second approach.
As far as we know, there is currently no other numerical method available for this problem.
We demonstrate the applicability and computational efficiency of the new method with several numerical tests, illustrating the convergence of the approximated eigenvalues to the exact ones and supporting the conjecture of spectral accuracy.

\bigskip

The paper is organized as follows.
In section \ref{s:model} we introduce the prototype model and the relevant solution operators and IG.
Section \ref{sec:reformulation2d} describes the reformulation of the IG in terms of the integrated state.
The reformulated IG is discretized in section \ref{s_numerical} and the resulting numerical method is applied to some test problems in section \ref{s_tests}.
In section \ref{s_velocity} we discuss the extension to models with structuring variables evolving with nontrivial velocities and show an example.
Finally, we provide some concluding remarks in section \ref{s_concluding}.
In appendix \ref{sec:1d}, for completeness, we apply the method to the case of one structuring variable and illustrate it with a test model.

A MATLAB implementation of the method is available at \url{http://cdlab.uniud.it/software}.

\section{The prototype model}
\label{s:model}

Let $x_{0},\bar x,y_{0},\bar y\in\mathbb{R}$ such that $x_{0}<\bar x$ and $y_{0}<\bar y$ and let $\overline{R} \coloneqq [x_0,\bar{x}]\times[y_0,\bar{y}]$.
We consider the scalar first-order linear hyperbolic PDE
\begin{equation}
\label{ruanPDE}
\partial_{t} u(t,x,y)+\partial_{x}u(t,x,y)+\partial_{y}u(t,x,y)=
- \mu(x,y)u(t,x,y),
\end{equation}
with boundary conditions
\begin{align}
u(t,x,y_{0})&=\iint_{\overline{R}}\alpha(x,\xi,\sigma)u(t,\xi,\sigma)\dd\xi\dd\sigma \eqqcolon K_{\alpha}(u(t,\cdot,\cdot))(x),\label{ruanBCa}\\
u(t,x_{0},y)&=\iint_{\overline{R}}\beta(y,\xi,\sigma)u(t,\xi,\sigma)\dd\xi\dd\sigma \eqqcolon K_{\beta}(u(t,\cdot,\cdot))(y),\label{ruanBCb}
\end{align}
where $u(t,x,y)$ is the density of the given population at time $t\geq0$ depending on the two structuring variables $x\in[x_{0},\bar x]$ and $y\in[y_{0},\bar y]$.

Following \cite[Assumption 2.1]{KangHuoRuan2021}, we assume that the model coefficients $\mu \colon \overline{R} \to \mathbb{R}$, $\alpha \colon [x_0,\bar{x}] \times \overline{R} \to \mathbb{R}$ and $\beta \colon [y_0,\bar{y}] \times \overline{R} \to \mathbb{R}$ are nonnegative $L^1$ functions, Lipschitz continuous on the interior of their domains, with $\alpha(x, \cdot, \cdot)$ and $\beta(y, \cdot, \cdot)$ dominated by $L^1$ functions in $x$ and $y$, respectively, and $\mu$ bounded from above and bounded away from $0$ (for a similar but more general approach, see \cite[section 4]{SinestrariWebb1987}).
Observe that $\alpha(x, \cdot, \cdot)$ and $\beta(y, \cdot, \cdot)$ are essentially bounded, so the operators $K_\alpha$ and $K_\beta$ map $L^1(\overline{R})$ to $L^1([x_0,\bar{x}])$ and $L^1([y_0,\bar{y}])$, respectively, and are bounded.

Under these assumptions, for every $u_0 \in L^1(\overline{R})$ the initial--boundary value problem defined by \eqref{ruanPDE}--\eqref{ruanBCb} and $u(0,\cdot,\cdot)= u_{0}$ admits a unique solution $u(t,\cdot,\cdot) \in L^1(\overline{R})$ for $t \geq 0$.
Moreover, the family of solution operators $\{T(t)\}_{t\geq 0}$, defined by $T(t)u_0 = u(t,\cdot,\cdot)$, forms a strongly continuous semigroup of bounded linear operators in the Banach space $L^1(\overline{R})$, see \cite[Theorem 3.1]{Webb1985c} or \cite[Theorem 2.3]{KangHuoRuan2021}.
In addition, $\{T(t)\}_{t\geq 0}$ is eventually compact \cite[section 4]{KangHuoRuan2021}.

The IG of the semigroup $\{T(t)\}_{t\geq 0}$ is the operator $A \colon D(A) \to L^1(\overline{R})$ defined by
\begin{equation}
\label{IGdef}
A\phi = \lim_{t \to 0^+} \frac{1}{t} ( T(t)\phi - \phi),
\end{equation}
where $D(A) \subset L^1(\overline{R})$ consists of the functions for which the limit exists.
Kang et al.\ \cite[Remark 2.3]{KangHuoRuan2021} prove that the operator $A$ satisfies
\begin{equation}
\label{IG-action}
\begin{aligned}
A\phi(x,y) &= - \partial_x \left[ \phi(x,y) + \partial_y \int_{x_0}^x \phi(a,y)\dd a \right] - \mu(x,y) \phi(x,y) \\
	&= - \partial_y \left[ \phi(x,y) + \partial_x \int_{y_0}^y \phi(x,a)\dd a \right] - \mu(x,y) \phi(x,y),
\end{aligned}
\end{equation}
for a.e.\ $x \in [x_0, \bar{x}]$, a.e.\ $y \in [y_0, \bar{y}]$, and every $\phi \in D(A)$, and that $D(A)$ satisfies the inclusion
\begin{align*}
D(A) \subset \biggl\{ &\phi \in L^1(\overline{R}) \;\bigg\vert \\
	& (x,y) \mapsto \int_{x_0}^x \phi(s,y) \dd s \text{ is absolutely continuous in $y$, for a.e. $x \in [x_0, \bar{x}]$}, \\
	& (x,y) \mapsto \left[ \phi(x,y) + \partial_y \int_{x_0}^x \phi(s,y) \dd s \right] \text{ is absolutely continuous in $x$,} \\
	& \hphantom{(x,y) \mapsto \left[ \phi(x,y) + \partial_y \int_{x_0}^x \phi(s,y) \dd s \right]{}} \text{ for a.e. $y \in [y_0, \bar{y}]$}, \\
	& \lim_{x \to x_0^+} \left[ \phi(x,y) + \partial_y \int_{x_0}^x \phi(s,y) \dd s \right] = K_\beta(\phi)(y) \text{ for a.e. $y \in [y_0, \bar{y}]$},  \\
	& (x,y) \mapsto \int_{y_0}^y \phi(x,s) \dd s \text{ is absolutely continuous in $x$, for a.e. $y \in [y_0, \bar{y}]$}, \\
	& (x,y) \mapsto \left[ \phi(x,y) + \partial_x \int_{y_0}^y \phi(x,s) \dd s \right] \text{ is absolutely continuous in $y$,} \\
	& \hphantom{(x,y) \mapsto \left[ \phi(x,y) + \partial_x \int_{y_0}^y \phi(x,s) \dd s \right]{}} \text{ for a.e. $x \in [x_0, \bar{x}]$}, \\
	& \lim_{y \to y_0^+} \left[ \phi(x,y) + \partial_x \int_{y_0}^y \phi(x,s) \dd s \right] = K_\alpha(\phi)(x) \text{ for a.e. $x \in [x_0, \bar{x}]$}, \\
	& \partial_x   \left[ \phi(x,y) + \partial_y \int_{x_0}^x \phi(s,y) \dd s \right] \in L^1(\overline{R}), \\
	& \partial_y   \left[ \phi(x,y) + \partial_x \int_{y_0}^y \phi(x,s) \dd s \right] \in L^1(\overline{R})
 \biggr\},
\end{align*}
while \cite[Remark 6.1]{Webb1985c} claims that equality holds.
If $\phi \in D(A)$ is sufficiently smooth, the action \eqref{IG-action} of the operator $A$ simplifies and can be expressed as
\begin{equation*}
A\phi = -(\partial_x \phi + \partial_y \phi + \mu \phi).
\end{equation*}

The spectrum of the IG determines the stability of the zero equilibrium.%
\footnote{Note that since $A$ is an operator on a real Banach space, in order to define and compute its spectrum the space and the operator need to be complexified.
For details see, e.g., \cite[section III.7]{DiekmannVanGilsVerduynLunelWalther1995}.}
More precisely, the latter is asymptotically stable if and only if the spectral abscissa of $A$ is negative and it is unstable if the spectral abscissa is positive (see \cite[Theorem 9.5]{ClementHeijmansAngenentVanDuijnDePagter1987} and \cite[Theorem VI.1.15]{EngelNagel2000}).

Observe that in the model defined by \eqref{ruanPDE}--\eqref{ruanBCb} the structuring variables evolve with the same velocity as time.
In sections \ref{sec:reformulation2d} and \ref{s_numerical} we present the integral reformulation and the discretization restricting to this special case, as in \cite{KangHuoRuan2021}.
However, they can be applied to the case of nontrivial velocities as described in section \ref{s_velocity}, where we also present an example.

\section{Equivalent formulation in a space of absolutely continuous functions}
\label{sec:reformulation2d}

To conveniently handle the boundary conditions in $D(A)$ from a numerical point of view, inspired by the approach of \cite{ScarabelBredaDiekmannGyllenbergVermiglio2021} in the case of one structuring variable, we argue in terms of the integrated state.
In particular, we define an isomorphism between $L^1(\overline{R})$ and a suitable space of functions via integration.
We then use this isomorphism and the semigroup $\{T(t)\}_{t \geq 0}$ to construct an appropriate semigroup acting on a space of functions with higher regularity.

With this goal in mind, we first recall the definition and some properties of absolute continuity in the sense of Carath\'eodory.
We refer the reader to \cite{Sremr2010} for further details.

For $(x,y)\in \overline{R}$, define $R(x,y)\coloneqq[x_0,x]\times[y_0,y]$; then $\overline{R} = R(\bar{x},\bar{y})$.
A function $v$ defined on $\overline{R}$ is \emph{absolutely continuous in the sense of Carathéodory} if and only if there exist $e_v\in\mathbb{R}$, $f_v\in L^1([x_0,\bar{x}])$, $g_v\in L^1([y_0,\bar{y}])$ and $h_v\in L^1(\overline{R})$ such that
\begin{equation*}
v(x,y) = e_v + \int_{x_0}^x f_v(a) \dd a + \int_{y_0}^y g_v(b) \dd b + \iint_{R(x,y)} h_v(a,b) \dd a \dd b.
\end{equation*}
Observe that the double integral in the last term is equal to the iterated integral on the variables $a$ and $b$ in any order, thanks to Fubini's theorem.
The space $AC(\overline{R})$ of absolutely continuous functions on $\overline{R}$ in the sense of Carathéodory is a Banach space when equipped with the norm $\lVert\cdot\rVert_{AC(\overline{R})}$ defined as
\begin{equation*}
\lVert v\rVert_{AC(\overline{R})} \coloneqq \lvert e_v\rvert + \lVert f_v\rVert_{L^1([x_0,\bar{x}])} + \lVert g_v\rVert_{L^1([y_0,\bar{y}])} + \lVert h_v\rVert_{L^1(\overline{R})}.
\end{equation*}
We consider a particular subspace of $AC(\overline{R})$, namely
\begin{equation*}
AC_0(\overline{R}) := \biggl\{ v \colon \overline{R} \to \mathbb{R} \;\bigg\vert\; v(x,y) = \iint_{R(x,y)} h_v(a,b) \dd a \dd b \ \text{ for some }h_v \in L^1(\overline{R})\biggr\}.
\end{equation*}
Observe that $AC_0(\overline{R})$ is a Banach space, being a closed subspace, and that $v(x_0,y)=v(x,y_0)=0$ for $v\in AC_0(\overline{R})$.
The operator $V \colon L^1(\overline{R}) \to AC_0(\overline{R})$ defined by
\begin{equation*}
V\phi(x,y) = \iint_{R(x,y)} \phi(a, b) \dd a \dd b
\end{equation*}
defines an isomorphism between $L^1(\overline{R})$ and $AC_0(\overline{R})$, with $V^{-1}\psi = \partial_x \partial_y \psi$ for all $\psi \in AC_0(\overline{R})$.
Observe that both $V$ and $V^{-1}$ are bounded ($\lVert V\rVert = \lVert V^{-1}\rVert = 1$).

Note that, given a solution $u$ of \eqref{ruanPDE}--\eqref{ruanBCb}, $\iint_{R(x,y)} u(t,a,b) \dd a \dd b$ represents the number of individuals whose structuring variables belong to $R(x,y)$ at time $t$.

Returning now to \eqref{ruanPDE}--\eqref{ruanBCb}, we define the family of operators $\{S(t)\}_{t \geq 0}$ on \linebreak $AC_0(\overline{R})$ as $S(t) \coloneqq VT(t) V^{-1}$.
Since $V$ and $V^{-1}$ are linear and bounded, recalling that the composition of a compact and a bounded operator (in either order) is compact, the operators $S(t)$ are in turn linear and bounded and form a family with the same properties as $\{T(t)\}_{t \geq 0}$, namely they form a strongly continuous and eventually compact semigroup on $AC_0(\overline{R})$.
Its IG is $B \colon D(B) \to AC_0(\overline{R})$, with $B\coloneqq VAV^{-1}$ and $D(B) = VD(A)$.

With the aim of using $B$ to study the stability properties of \eqref{ruanPDE}--\eqref{ruanBCb}, it is important to understand the relation between the spectra of $A$ and $B$.
By \cite[Corollary V.3.2(i)]{EngelNagel2000}, since the semigroups $\{T(t)\}_{t \geq 0}$ and $\{S(t)\}_{t \geq 0}$ are eventually compact, the spectra of both $A$ and $B$ are at most countable and consist only of eigenvalues (of finite algebraic multiplicity).
By \cite[Proposition 4.1]{BredaMasetVermiglio2012} $A$ and $B$ have the same nonzero eigenvalues (with the same multiplicities) and corresponding eigenvectors.
Moreover, observe that $A$ is injective if and only if $B$ is, so $0$ is either in both spectra or in none.
Therefore, the spectra of $A$ and $B$ coincide.
From \cite[Proposition 3.1]{KangHuoRuan2021}, noting that in that proof $\chi = A\phi + \mu \phi$, for $\psi \in D(B)$ we can write
\begin{equation}
\label{B-action}
\begin{aligned}
B\psi(x,y) = &- \partial_x \psi(x,y) - \partial_y \psi(x,y) \\
&+ \int_{x_0}^x K_\alpha(\partial_y\partial_x \psi)(a) \dd a + \int_{y_0}^y K_\beta(\partial_y\partial_x \psi)(b) \dd b \\
&- \iint_{R(x,y)} \mu(a,b) \partial_y\partial_x \psi (a,b) \dd a \dd b.
\end{aligned}
\end{equation}

\section{Pseudospectral discretization of the IG}
\label{s_numerical}

In this section, we use pseudospectral methods with a tensorial approach to obtain a finite-dimensional approximation of the operator $B$, whose spectrum can be used to determine the stability properties of the system.

For $n$ and $m$ positive integers, let $\Pi_{n,m}^0$ be the space of bivariate polynomials on $\overline{R}$ of degree at most $n$ in the first variable and at most $m$ in the second variable, taking value $0$ at $x=x_0$ and $y=y_0$.
These conditions are motivated by the fact that we use polynomials in $\Pi_{n,m}^0$ as approximations of functions in $AC_0(\overline{R})$.
Let $\Theta_X = \{x_1,\ldots,x_n\}$ be a mesh of $n$ points in $(x_{0},\bar x]$, with $x_{0}<x_{1}<\dots<x_{n}=\bar x$, and let $\Theta_Y = \{y_1,\ldots,y_m\}$ be a mesh of $m$ points in $(y_{0},\bar y]$, with $y_{0}<y_{1}<\dots<y_{m}=\bar y$.
We approximate a function $\psi \in AC_0(\overline{R})$ by a vector $\Psi \in \mathbb{R}^{nm}$ according to
\begin{equation*}
\psi(x_{i},y_{j})=\Psi_{i,j},\quad i=1,\ldots,n,\;j=1,\ldots,m,
\end{equation*}
where the components of $\Psi$ are ordered according to the lexicographic order of the double indices $(i,j)$.

Given $\Psi \in \mathbb{R}^{nm}$, let $\psi_{n,m}\in\Pi_{n,m}^0$ be the polynomial interpolating $\Psi$ on $\Theta_X \times \Theta_Y$:
\begin{align*}
\psi_{n,m}(x_{i},y_{j}) = \Psi_{i,j}, \quad i=1,\ldots,n,\;j=1,\ldots,m.
\end{align*}
The finite-dimensional approximation of the operator $B$ is then $B_{n,m} \colon \mathbb{R}^{nm} \to \mathbb{R}^{nm}$ defined as
\begin{equation}\label{Anm}
[B_{n,m}\Psi]_{i,j}\coloneqq (B\psi_{n,m})(x_{i},y_{j}), \quad i=1,\ldots,n,\;j=1,\ldots,m.
\end{equation}

We can write more explicitly the entries of the matrix $B_{n,m}$ by using the bivariate Lagrange representation of $\psi_{n,m}$, together with the explicit action of the operator $B$ defined in \eqref{B-action}.
Let $\{\ell_{X,i}\}_{i=0,\ldots,n}$ and $\{\ell_{Y,j}\}_{j=0,\ldots,m}$ be the Lagrange bases of polynomials relevant to $\{x_0\} \cup \Theta_X$ and $\{ y_0\} \cup \Theta_Y$, i.e.,
\begin{align*}
\ell_{X,i}(x) = \prod_{\substack{k=0\\k\neq i}}^n \frac{x-x_k}{x_i-x_k}, \quad 
\ell_{Y,j}(y) = \prod_{\substack{k=0\\k\neq j}}^m \frac{y-y_k}{y_j-y_k}.
\end{align*}
The polynomial $\psi_{n,m}$ can be written as
\begin{equation*}
\psi_{n,m}(x,y)=\sum_{i=1}^{n}\sum_{j=1}^{m}\ell_{X,i}(x)\ell_{Y,j}(y)\Psi_{i,j}, \quad (x,y)\in \overline{R}.
\end{equation*}
Note that indeed $\psi_{n,m}(x,y)=0$ for $x=x_0$ or $y=y_0$.
Using \eqref{B-action} and \eqref{Anm}, we get
\begin{equation*}
\begin{aligned}
[B_{n,m}\Psi]_{k,l} = 
& - \sum_{i=1}^{n} \ell_{X,i}'(x_k)\Psi_{i,l} - \sum_{j=1}^{m} \ell_{Y,j}'(y_l) \Psi_{k,j} \\
&+ \int_{y_0}^{y_l} K_{\beta}\biggl(\sum_{i=1}^{n}\sum_{j=1}^{m}\ell_{X,i}'\ell_{Y,j}' \Psi_{i,j}\biggr)(b) \dd b \\
&+ \int_{x_0}^{x_k} K_{\alpha}\biggl(\sum_{i=1}^{n}\sum_{j=1}^{m}\ell_{X,i}'\ell_{Y,j}'\Psi_{i,j}\biggr)(a) \dd a \\
&- \iint_{R(x_k,y_l)} \mu(a,b) \sum_{i=1}^{n}\sum_{j=1}^{m}\ell_{X,i}'(a) \ell_{Y,j}'(b) \Psi_{i,j} \dd a \dd b.
\end{aligned}
\end{equation*}
Using the linearity of $K_\alpha$ and $K_\beta$, it is easy to characterize the entries of the matrix $B_{n,m}$. 
Let $D_X\in \mathbb{R}^{n\times n}$ and $D_Y\in \mathbb{R}^{m\times m}$ be defined as
\begin{alignat*}{2}
[D_X]_{i,j} &= \ell_{X,j}'(x_i), &\quad i,j&=1,\ldots,n \\
[D_Y]_{i,j} &= \ell_{Y,j}'(y_i), &\quad i,j&=1,\ldots,m.
\end{alignat*}
In other words, $D_X$ and $D_Y$ are the part of the differentiation matrices associated with $\{x_0\} \cup \Theta_X$ and $\{ y_0 \} \cup \Theta_Y$, respectively, deleting the first row and the first column. The bivariate differentiation matrices in $x$ and $y$ are
\begin{equation*}
\mathbf{D}_X=D_X \otimes I_m, \quad \mathbf{D}_Y=I_n \otimes D_Y,
\end{equation*}
where $\otimes$ denotes the Kronecker product.
We can then write
\begin{equation*}
B_{n,m} = - \mathbf{D}_X - \mathbf{D}_Y + \mathbf{A} + \mathbf{B} - \mathbf{M},
\end{equation*}
where $\mathbf{A},\mathbf{B},\mathbf{M} \in \mathbb{R}^{nm \times nm}$ are defined by
\begin{align}
\label{A}
\mathbf{A}_{(k,l),(i,j)} &= \int_{x_0}^{x_k} K_\alpha(\ell_{X,i}'\ell_{Y,j}')(a)\dd a, \\
\label{B}
\mathbf{B}_{(k,l),(i,j)} &= \int_{y_0}^{y_l} K_\beta(\ell_{X,i}'\ell_{Y,j}')(b)\dd b, \\
\label{M}
\mathbf{M}_{(k,l),(i,j)} &= \iint_{R(x_k,y_l)} \mu(a,b) \ell_{X,i}'(a) \ell_{Y,j}'(b) \dd a \dd b,
\end{align}
for $k,i=1,\ldots,n$ and $l,j=1,\ldots,m$. 
Note that, if $\mu$ is constant, the matrix $\mathbf{M}$ is diagonal with diagonal entries equal to $\mu$.

\bigskip

We finally note that, although the matrix $B_{n,m}$ is defined for any set of nodes, the choice of the latter is critical to ensure the convergence of the interpolating polynomials and, in turn, of the elements of the spectrum. 
In the following numerical experiments, we choose the Chebyshev extremal points in each interval $[x_0,\bar x]$ and $[y_0,\bar y]$.

In the univariate case, these nodes guarantee that the convergence rate of the interpolating polynomial of degree $n$ is $O(n^{-k})$ if the interpolated function is $C^k$ \cite[Theorem 7.2]{Trefethen2013}, which implies that the order of convergence is infinite if the function is smooth.
Moreover, the convergence rate is $O(c^n)$ for some $0<c<1$ if the function is analytic \cite[Theorem 8.2]{Trefethen2013}.
The two latter properties are often known as \emph{spectral accuracy}, see \cite[chapter 4]{Trefethen2000} and \cite[chapter 2]{Boyd2001}.
Furthermore, observe that the relevant differentiation matrices can be computed explicitly \cite{Trefethen2000}.

The classic result on the interpolation error being bounded by means of the best uniform approximation error and the Lebesgue constant holds also in the bivariate case.
A multidimensional version of Jackson's theorem on the best uniform approximation error holds as well \cite[Theorem 4.8]{Schultz1969}.
Moreover, it is easy to verify that the Lebesgue constant for the tensor-product Chebyshev extremal nodes in $\overline{R}$ is the product of the univariate Lebesgue constants in $[x_0, \bar x]$ and $[y_0, \bar y]$, hence it is $O(\log n \log m)$.
The tensor-product Chebyshev interpolation is thus near-optimal also in the bivariate case.

Although a proof of convergence for the method is out of the scope of this paper, we show that the order of convergence observed numerically for the approximated eigenvalues and eigenvectors is consistent with the well-established order of convergence of polynomial interpolation.

\bigskip

For implementation purposes, we observe that in general the integrals defining $\mathbf{A}$, $\mathbf{B}$ and $\mathbf{M}$ in \eqref{A}--\eqref{M} cannot be computed exactly. 
To approximate the integrals on $\overline{R}$ we use the Clenshaw--Curtis cubature formula \cite{ClenshawCurtis1960}, which is based on Chebyshev extremal points and is spectrally accurate \cite{Trefethen2000,Trefethen2008}.

To compute the entries of $\mathbf{A}$ and $\mathbf{B}$, given a function $f$ defined on $[x_0,\bar{x}]$ and $F \in \mathbb{R}^n$ such that $F_i = f(x_i)$, we consider the approximation
\begin{equation*}
\int_{x_0}^{x_k} f(a) \dd a \approx [D_X^{-1} F]_k, \quad k=1,\ldots,n,
\end{equation*}
and similarly for functions defined on $[y_0,\bar{y}]$ (see, e.g., \cite{DiekmannScarabelVermiglio2020}).
This approximation can be extended to the double integrals involved in the entries of the matrix $\mathbf{M}$ for nonconstant $\mu$. More precisely, given a function $\psi \in AC_0(\overline{R})$ and a vector $\Psi \in \mathbb{R}^{nm}$ such that $\Psi_{k,l} = \psi(x_k,y_l)$, each integral $\int_{x_0}^{x_k} \int_{y_0}^{y_l} \psi(a,b)\dd a \dd b$ can be approximated by the corresponding $(k,l)$-th entry of the vector $\mathbf{D}_X^{-1}\mathbf{D}_Y^{-1}\Psi$.

\section{Numerical experiments}
\label{s_tests}

In this section, we present several numerical experiments to investigate how the spectrum of the finite-dimensional operator $B_{n,m}$ approximates the spectrum of $B$, and in turn of $A$, of each problem at hand. For this purpose, we select several parameter sets for which eigenvalues and eigenfunctions of $A$ can be expressed analytically, and we study the convergence of the approximated eigenvalues of $B_{n,m}$ to the analytic ones.
As for the eigenfunctions, we stress that, since $B_{n,m}$ represents an approximation of the operator $B$, an eigenvector $\Psi$ of $B_{n,m}$ provides an approximation $\psi_{n,m}$ of $V\phi$, where $\phi$ is an eigenfunction of $A$; an approximation of $\phi$ is thus given by $\partial_x\partial_y \psi_{n,m}$.

For each example we study the behavior for increasing $n=m$ of the absolute error $\epsilon_\lambda$ on the known eigenvalue $\lambda$ and of the absolute error $\epsilon_\phi$ in $L^1$ norm on the known eigenfunction $\phi$, computed via Clenshaw--Curtis cubature.

In all examples we choose
\begin{equation*}
\alpha(x, \xi,\sigma)= \alpha_1(x)\gamma(\xi, \sigma), \quad \beta(y, \xi, \sigma)=\beta_1(y)\gamma(\xi, \sigma)
\end{equation*}
in the boundary conditions \eqref{ruanBCa}--\eqref{ruanBCb}, in order to simplify finding an explicit eigenfunction.

We remark that the parameters are chosen in order to have an analytically known eigenfunction with certain smoothness properties, without regard to any specific biological interpretation.

To compute the spectrum of $B_{n,m}$ we use standard methods (namely MATLAB's \texttt{eig} function).
Note that the approximated spectrum may contain spurious eigenvalues (e.g., when $B$ has fewer eigenvalues than the dimension of $B_{n,m}$); however, in our examples we only examine specific eigenvalues, so that the possible spurious ones do not affect our analysis.

\subsection{Analytic eigenfunctions}

\begin{table}
\centering
{\renewcommand*\arraystretch{1.5}
\rowcolors{2}{gray!10}{white}
\begin{tabular}{cccccc}
\rowcolor{gray!20}
  & Ex.\ 1.1 & Ex.\ 1.2 & Ex.\ 1.3 & Ex.\ 1.4\\ 
 \specialrule{\lightrulewidth}{0pt}{0pt}
 $\overline{R}$  & $[0,1]^2$  & $\left[\frac{\pi}{6}, \frac{\pi}{2}\right]\times\left[\frac{\pi}{6},\frac{\pi}{4}\right]$ & $[0, 2]\times[{-1},1]$ & $[0, 2]\times[0,1]$ \\
 $\alpha_1(x)$ & $1$ & $\cos(x-\frac{\pi}{6})$  & $e^{x+1}$  & $e^{-x^2}$ \\
 $\beta_1(y)$  & $1$ &$\cos(\frac{\pi}{6}-y)$  & $e^{-y}$  & $e^{y}$ \\
 $\gamma(\xi, \sigma)$ & $1$ & $\left(\frac{\sqrt 2-\sqrt6+2}{4}\right)^{-1}$  & $\frac{1}{4}e^{-\xi+\sigma}$  & $\left(\iint_{\overline{R}} e^{-a^2+b}\dd b\dd a\right)^{-1}$ \\
 $\mu(x, y)$  & $1$ & $1$ & $1$ &  $2x+1$  \\
 \specialrule{\lightrulewidth}{0pt}{0pt}
 $\lambda$ & $-1$ & $-1$ & $-1$  & $-2$ \\
 $\phi(x, y)$  & $1$ &$\cos(x-y)$  & $e^{x-y}$  & $e^{-x^2+y}$
\end{tabular}
}
\caption{Parameters and resulting eigenvalue and eigenfunction in the analytic cases.}
\label{tab1}
\end{table}

We consider a first group of examples presenting an analytic eigenfunction.
The choices of the parameters and the resulting eigenvalue and eigenfunction are listed in Table \ref{tab1}.
Starting from Example 1.1, where all parameters are constant, we gradually introduce nonconstant coefficients: $\alpha_1$ and $\beta_1$ in Example 1.2, $\gamma$ in Example 1.3 and $\mu$ in Example 1.4.

\begin{figure}
\centering
\includegraphics[trim=28.668pt 0 7.353pt 0]{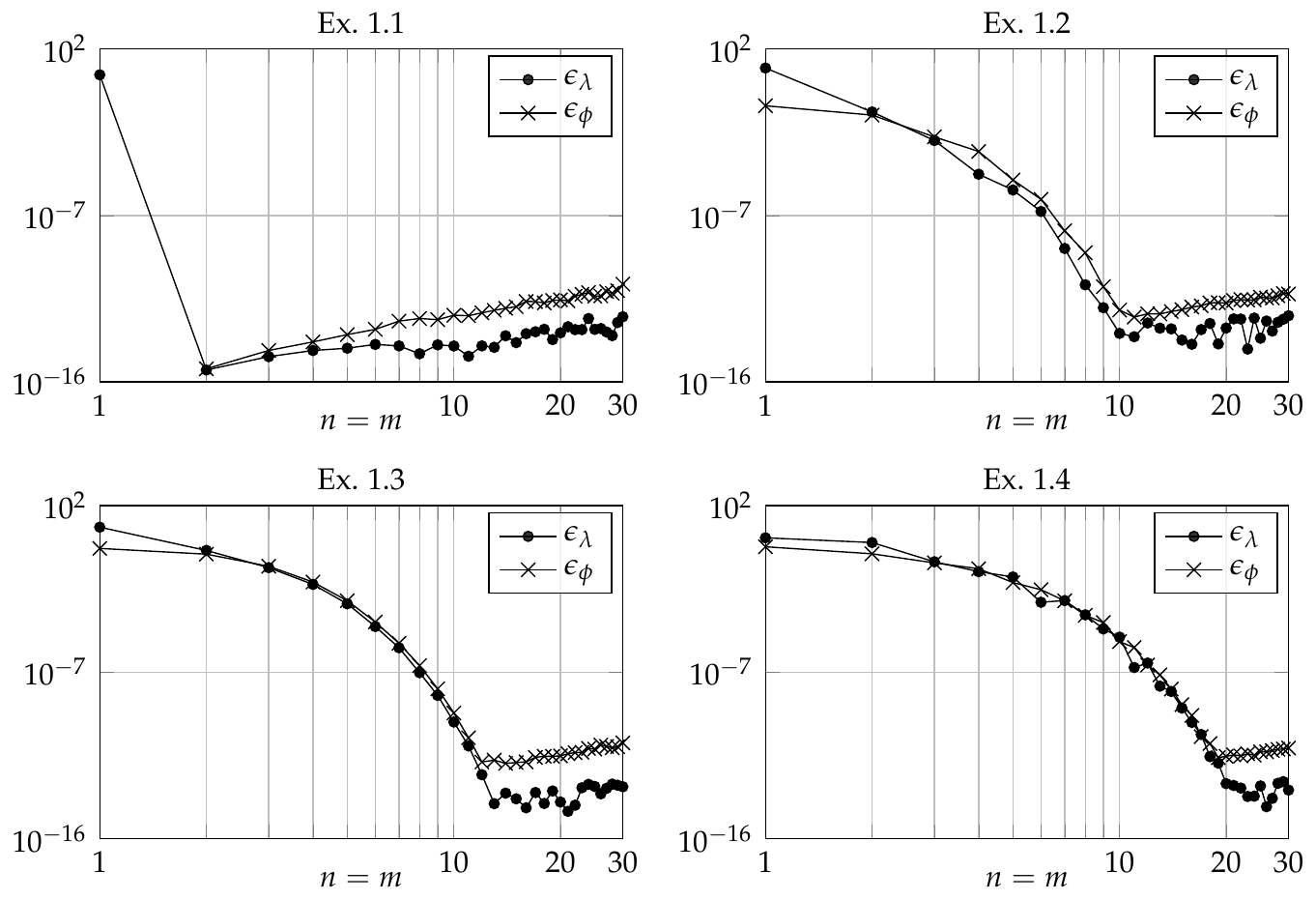}
\caption{Errors $\epsilon_\lambda$ and $\epsilon_\phi$ for Examples 1.1, 1.2, 1.3 and 1.4 defined in Table~\ref{tab1}.
Observe that for Example 1.1 with $n=m=1$ the error $\epsilon_\phi$ is exactly $0$, hence it is not represented in the logarithmic scale.}
\label{fig:01-02-03-04}
\end{figure}

Considering Example 1.1, Figure \ref{fig:01-02-03-04} shows that the errors reach the machine precision already for $n=m=2$.
With $n=m=1$ the error $\epsilon_\phi$ on the eigenfunction is exactly equal to $0$, which may be explained by the fact that constant functions are interpolated exactly already by polynomials of degree $0$.
As $n=m$ increases, the errors increase, possibly due to numerical instability.

Considering now Examples 1.2, 1.3 and 1.4, Figure \ref{fig:01-02-03-04} shows that both $\epsilon_\lambda$ and $\epsilon_\phi$ decay with infinite order.
Observe that in Example 1.4 more nodes are needed to reach the error barrier than in Examples 1.2 and 1.3, probably due to the approximation of the integrals involving the nonconstant $\mu$.

\subsection{Nonsmooth eigenfunctions}

\begin{table}
\centering
{\renewcommand*\arraystretch{1.5}
\rowcolors{2}{gray!10}{white}
\begin{tabular}{cccccc}
\rowcolor{gray!20}
 & Ex.\ 2.1 ($C^{2}$) & Ex.\ 2.2 ($C^{1}$) & Ex.\ 2.3 ($C^{0}$) & Ex.\ 2.4 (discontinuous) \\ 
 \specialrule{\lightrulewidth}{0pt}{0pt}
 $\overline{R}$ &$[0, 1]\times[0,2]$  & $[0, 1]\times[0,2]$ & $[0, 1]\times[0,2]$ & $[0,1]\times[0,2]$ \\
 $\alpha_1(x)$ & $x^2\lvert x\rvert$  & $-x\lvert x\rvert$  & $\lvert x\rvert$  & $\mathbf{1}_{[0, +\infty[}(x)$ \\
 $\beta_1(y)$  & $y^2\lvert y\rvert$  & $y\lvert y\rvert$  & $\lvert y\rvert$  & $\mathbf{1}_{]-\infty, 0]}(y)$ \\
 $\gamma(\xi, \sigma)$ & $\frac{5}{8}$  & $\frac{6}{7}$ &  $\frac{3}{4}$   & $2$ \\
 $\mu(x, y)$  & $1$ & $1$ & $1$ &  $1$  \\
 \specialrule{\lightrulewidth}{0pt}{0pt}
 $\lambda$ & $-1$ & $-1$ & $-1$  & $-1$ \\
 $\phi(x, y)$ & $(x-y)^2\lvert x-y\rvert$  &$(x-y)\lvert x-y\rvert$  & $\lvert x-y\rvert$  & $\mathbf{1}_{[0, +\infty[}(x-y)$
\end{tabular}
}
\caption{Parameters and resulting eigenvalue and eigenfunction in the nonsmooth case.
$\mathbf{1}_{A}$ is the indicator function of $A\subset \mathbb{R}$.
The regularity of the eigenfunction is shown in parentheses in the titles.}
\label{tab2}
\end{table}

For the second group of examples, we consider eigenfunctions which are not smooth.
The choices of the parameters and the resulting eigenvalue and eigenfunction are listed in Table \ref{tab2}.
Observe that the eigenfunctions have the same regularity as the coefficients $\alpha_1$ and $\beta_1$, namely $C^{2}$ for Example 2.1, $C^1$ for Example 2.2 and $C^{0}$ for Example 2.3.

We consider also Example 2.4 with discontinuous coefficients and eigenfunction.
In \cite{KangHuoRuan2021} it is assumed that $\alpha$ and $\beta$ are Lipschitz continuous in order to prove that the semigroup $\{T(t)\}_{t\geq0}$ is eventually compact, but actually this hypothesis is used only in the interior of the domains: thus, Example 2.4 still falls inside the scope of \cite{KangHuoRuan2021}.
Observe also that the function $\mathbf{1}_{[0, +\infty[}(x-y)$, even if not continuous, is in $D(A)$.

\begin{figure}
\centering
\includegraphics[trim=25.016pt 0 7.353pt 0]{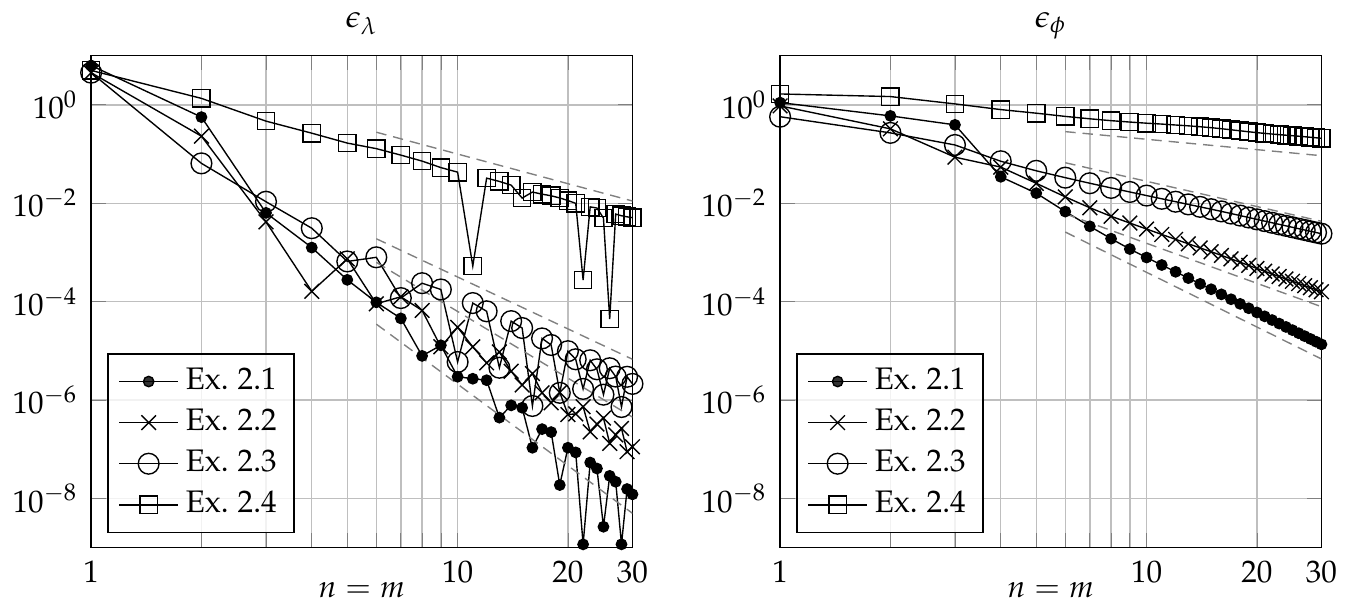}
\caption{Errors $\epsilon_\lambda$ (left) and $\epsilon_\phi$ (right) for Examples 2.1, 2.2, 2.3 and 2.4 defined in Table~\ref{tab2}.
The slopes of the dashed gray lines are $-2$, $-3.5$, $-4.5$, $-5.5$ (left) and $-0.7$, $-1.7$, $-2.7$, $-3.7$ (right), included to ease the interpretation of the plots.}
\label{fig:10-06-07-12-err}
\end{figure}

Figure \ref{fig:10-06-07-12-err} suggests that the errors decay with finite order and that these orders increase with the regularity of the eigenfunction.
In particular, focusing on Examples 2.1--2.3, we can observe that for both errors a loss of one order of differentiability of the eigenfunction seems to correspond to a loss of about one order of convergence (cfr. the dashed reference lines in Figure \ref{fig:10-06-07-12-err}).
The convergence of $\epsilon_\lambda$ seems to be almost two orders faster than that of $\epsilon_\phi$.
To possibly explain this difference, recall that we are actually collocating the eigenvalue problem for $B$, which means that the eigenvalues are the same as $A$, but the eigenfunctions correspond to integrals of the eigenfunctions of $A$, so the comparison between the eigenfunctions involves differentiating the computed ones.

As for Example 2.4, it appears that losing continuity of the eigenfunction itself changes the order of convergence differently for $\epsilon_\lambda$ and $\epsilon_\phi$.

\section{Structuring variables with nontrivial velocity}
\label{s_velocity}

We have illustrated the method for systems in which both physiological variables evolve at the same velocity as time.
This should not be seen as too restrictive, as systems with more general velocity terms \cite{BekkalBrikciClairambaultRibbaPerthame2008,DysonVillellaBressanWebb2000,Howard2003,MagalRuan2008,MetzDiekmann1986} can in some cases be reduced to \eqref{ruanPDE}--\eqref{ruanBCb} after a suitable scaling of variables, so that similar theoretical results on the stability of the zero solution hold \cite{Webb1985c}.
In practice, however, the computation of the change of variables, which in general is defined by the solution of an ODE system, may be expensive, although necessary when the individual parameters (e.g., birth and mortality rates) depend on the original (unscaled) variables. In this case, directly approximating the original problem with nontrivial velocities may be convenient from a computational point of view, as observed in \cite{ScarabelBredaDiekmannGyllenbergVermiglio2021}.

In fact, the transformation via integration can be easily carried out for problems of the form
\begin{gather*}
\partial_t u(t,x,y) + \partial_x (g_X(x) u(t,x,y)) + \partial_y (g_Y(y) u(t,x,y)) = - \mu(x,y) u(t,x,y), \\
\begin{aligned}
g_Y(y_0) u(t,x,y_0) &= K_\alpha(u(t,\cdot,\cdot))(x),\\
g_X(x_0) u(t,x_0,y) &= K_\beta(u(t,\cdot,\cdot))(y),
\end{aligned}
\end{gather*}
where the positive functions $g_X(x)$ and $g_Y(y)$ describe the rates of change of $x$ and $y$ in time.
In this case, it is straightforward to verify that, given the IG $A$, the operator $B=VAV^{-1}$ admits the representation
\begin{equation*}
\begin{aligned}
B\psi(x,y) = &- g_X(x)\partial_x \psi(x,y) - g_Y(y)\partial_y \psi(x,y) \\
&+ \int_{x_0}^x K_\alpha(\partial_y\partial_x \psi)(a) \dd a + \int_{y_0}^y K_\beta(\partial_y\partial_x \psi)(b) \dd b \\
&- \iint_{R(x,y)} \mu(a,b) \partial_y\partial_x \psi(a,b) \dd a \dd b,
\end{aligned}
\end{equation*}
which is approximated by a matrix of the form
\begin{equation*}
    B_{n,m} = - \mathbf{G}_X \mathbf{D}_X - \mathbf{G}_Y \mathbf{D}_Y + \mathbf{A} + \mathbf{B} - \mathbf{M},
\end{equation*}
where $\mathbf{A}$, $\mathbf{B}$ and $\mathbf{M}$ are the matrices defined in section \ref{s_numerical} and $\mathbf{G}_X=G_X \otimes I_m$, $\mathbf{G}_Y=I_n \otimes G_Y$, with $G_X$ and $G_Y$ diagonal matrices defined by $[G_{X}]_{i,i} = g_X(x_i)$, $i=1,\dots,n$, and $[G_{Y}]_{i,i} = g_Y(y_j)$, $j=1,\dots,m$.

\begin{figure}
\centering
\includegraphics[trim=28.668pt 0 7.353pt 0]{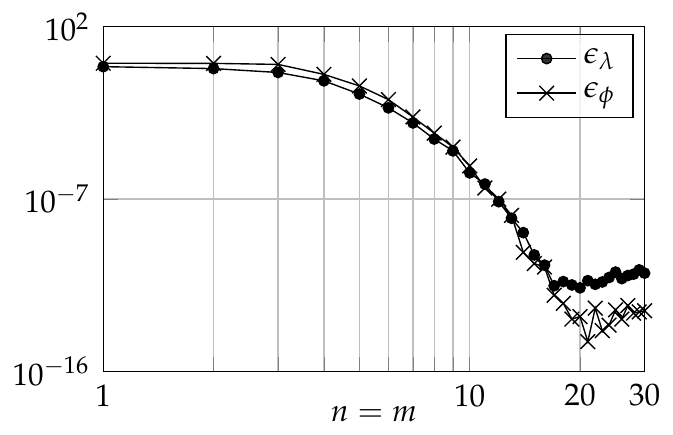}
\caption{Errors $\epsilon_\lambda$ and $\epsilon_\phi$ for the example of section \ref{s_velocity}.}
\label{fig:speed}
\end{figure}

As an example let us consider
\begin{gather*}
\overline{R}\coloneqq \left[\frac{1}{2}, \frac{3}{2}\right]\times \left[\frac{1}{2}, 2\right], \quad
g_X(x,y) \coloneqq x, \quad
g_Y(x,y) \coloneqq \frac{y^2}{2}, \\
\mu(x,y) \coloneqq y^3-2x^2-y+4, \\
\alpha(x, \xi, \sigma) \coloneqq \frac{e^{x^2-\frac{1}{4}}}{8c}, \quad
\beta(y, \xi, \sigma) \coloneqq \frac{e^{-y^2+\frac{1}{4}}}{2c}, \quad
c=\frac{1}{4} \pi \iint_{\overline{R}} e^{a^2-b^2} \dd a \dd b,
\end{gather*}
for which $\lambda=-5$ is an eigenvalue of the corresponding IG with eigenfunction $\phi(x, y) = e^{x^2-y^2}$. 
We can observe in Figure \ref{fig:speed} that the errors computed by our method decay with infinite order even in this case.

\section{Concluding remarks}
\label{s_concluding}

In this paper we proposed a numerical technique to analyze the stability of the zero solution of linear population models with two structuring variables, of the type considered in \cite{KangHuoRuan2021,Webb1985c}.

Extensive numerical tests illustrate the convergence of the eigenvalues of the finite-dimensional approximation to the true eigenvalues of the IG. 
The numerical tests support the conjecture that the order of convergence of the approximation depends on the regularity of the eigenfunction. A rigorous theoretical proof of the convergence of the approximation is left to future work.

Stability analysis requires not only that the eigenvalues of the IG are approximated accurately, but also that no spurious eigenvalues are to the right of the true spectral abscissa.
In fact, in our examples we observe that the approximation of the known eigenvalue is the numerical spectral abscissa, suggesting that the method can be effectively used to study the stability.

\bigskip
Structured population models can also be formulated as renewal equations for the population birth rate (or ``recruitment function'') \cite{BarrilCalsinaDiekmannFarkas2021,CalsinaDiekmannFarkas2016,DiekmannGyllenbergMetzNakaokaDeRoos2010}.
The renewal equation formulation is particularly convenient from the theoretical point of view as one can exploit the principle of linearized stability for nonlinear equations, which can not be proved in general for the PDE formulation \cite{BarrilCalsinaDiekmannFarkas2021}.
Results on the asymptotic behavior of solutions have also been recently proved, under special assumptions, for renewal equations defined on a space of measures, which makes it possible to consider a wider set of solutions compared to PDEs \cite{FrancoDiekmannGyllenberg2021,FrancoGyllenbergDiekmann2021}.

It would be interesting to apply pseudospectral methods in the framework of renewal equations (admitting a state space of multivalued functions or even measures), by extending the techniques developed for scalar renewal equations \cite{BredaDiekmannMasetVermiglio2013,BredaLiessi2018,ScarabelDiekmannVermiglio2021}.
However, when the structuring variables evolve with nontrivial velocity, the renewal equation formulation requires to explicitly invert the age-structure relation defined implicitly by an ODE system, which suffers from the computational challenges highlighted in section \ref{s_velocity}.
Hence, as explained therein, directly tackling the PDE formulation may be computationally convenient, as it bypasses the solution of the ODE system.

\bigskip
In this paper we restricted to structuring variables in bounded intervals, as this allows to exploit the highly desirable convergence properties of polynomial interpolation on Chebyshev nodes.
However, unbounded domains are common in the modeling literature (e.g., \cite{BekkalBrikciClairambaultRibbaPerthame2008,DysonVillellaBressanWebb2000,Inaba2017,SinestrariWebb1987}), for instance when it is not easy to determine a suitable upper bound for a physiological variable \emph{a priori}, or when the processes are naturally described by probability distributions with unbounded support (e.g., exponential or Gamma distributions).

In order to numerically treat these problems, truncating the domain can be feasible sometimes, but the accuracy of the approximation would depend on both the size of the truncated domain and the number of nodes in the domain. The latter usually becomes very large because the choice of Chebyshev nodes does not exploit the specific characteristics of the solutions, which usually belong to exponentially weighted spaces \cite{BarrilCalsinaDiekmannFarkas2021}.
For this reason, using exponentially weighted interpolation and Laguerre-type nodes has proved successful and more efficient than domain truncation in the case of delay equations \cite{GyllenbergScarabelVermiglio2018,ScarabelVermiglio}.
It would be interesting to apply similar techniques \cite{Boyd2001} to structured models in the PDE formulation with one or even two structuring variables, although the latter brings in additional complications due to the necessity to rely on multivariate interpolation.

\appendix

\section{Models with one structuring variable}
\label{sec:1d}

As recalled in the introduction, \cite{BredaCusulinIannelliMasetVermiglio2007} already provides a pseudospectral method, based on a different approach, to approximate the IG of models with one structuring variable.
In this appendix, for completeness, we adapt our approach to this case, providing also an example.

We consider the scalar first-order linear hyperbolic PDE
\begin{equation}
\label{1d-PDE-u}
\partial_{t} u(t,x)+\partial_{x} u(t,x) =
- \mu(x) u(t,x),
\end{equation}
with boundary condition
\begin{equation}
\label{1d-BC-u}
u(t,x_{0})= \int_{x_0}^{\bar{x}}\beta(\sigma)u(t,\sigma)\dd\sigma \eqqcolon K_\beta (u(t,\cdot)).
\end{equation}
As references on single structure models, see \cite{Iannelli1995,Inaba2017,Webb1985}; in particular, see \cite[section 1.2]{Inaba2017} for what concerns this appendix.

If $\mu, \beta \in L^1([x_0, \bar{x}])$ are nonnegative and bounded, for every $u_0 \in L^1([x_0, \bar{x}])$ the initial--boundary value problem defined by \eqref{1d-PDE-u}--\eqref{1d-BC-u} and $u(0,\cdot)= u_{0}$ admits a unique solution $u(t,\cdot) \in L^1([x_0,\bar{x}])$ for $t \geq 0$.
Moreover, the family of solution operators $\{T(t)\}_{t\geq 0}$, defined by $T(t)u_0 = u(t,\cdot)$, forms a strongly continuous and eventually compact semigroup of bounded linear operators in the Banach space $L^1([x_0,\bar{x}])$.
Its IG $A \colon D(A) \to L^1([x_0, \bar{x}])$, defined as in \eqref{IGdef}, can be expressed as
\begin{equation*}
A\phi = - \phi' - \mu\phi,
\end{equation*}
and its domain $D(A)$ can be characterized as
\begin{equation*}
D(A) = \{ \phi \in AC([x_0, \bar{x}]) \;\vert\; \phi(x_0) = K_\beta (\phi) \}.
\end{equation*}

Let us equip $AC([x_0, \bar{x}])$ with the norm $\lVert\cdot\rVert_{AC([x_0,\bar{x}])}$ defined as $\lVert f\rVert_{AC([x_0,\bar{x}])} \coloneqq \lvert f(x_0)\rvert + \lVert f'\rVert_{L^1([x_0,\bar{x}])}$.
Let $AC_0([x_0,\bar{x}])$ be the subspace of $AC([x_0, \bar{x}])$ of functions that are null at $x_0$, which is a Banach space, being a closed subspace.
The operator $V \colon L^1([x_0,\bar{x}]) \to AC_0([x_0,\bar{x}])$ defined by $V\phi(x) = \int_{x_0}^x \phi(\sigma) \dd \sigma$ defines an isomorphism between $L^1([x_0,\bar{x}])$ and $AC_0([x_0,\bar{x}])$, with $V^{-1}\psi = \psi'$ for all $\psi \in AC_0([x_0,\bar{x}])$.
Observe that both $V$ and $V^{-1}$ are bounded ($\lVert V\rVert = \lVert V^{-1}\rVert = 1$).

We define the family of operators $\{S(t)\}_{t \geq 0}$ on $AC_0([x_0,\bar{x}])$ as $S(t) \coloneqq VT(t) V^{-1}$.
As in section \ref{sec:reformulation2d}, we observe that they form a strongly continuous and eventually compact semigroup on $AC_0([x_0,\bar{x}])$ with IG $B \colon D(B) \to AC_0([x_0,\bar{x}])$, with $B\coloneqq VAV^{-1}$ and $D(B) = VD(A)$.
We can also derive the following expression for $B$, given $\psi \in D(B)$:
\begin{equation*}
B\psi = - \psi' + K_\beta (\psi') - V(\mu\psi').
 \end{equation*}
As in section \ref{sec:reformulation2d}, we can conclude that $A$ and $B$ have the same spectrum, at most countable and consisting only of eigenvalues (of finite algebraic multiplicity).

\bigskip

As an example, let us choose $[x_0, \bar x]=[0,2]$, $\mu(x)\equiv 1$ and $\beta(x)\coloneqq e^{-x}$.
It can be shown that the only real eigenvalue of the corresponding IG is the unique real solution of the equation
\begin{equation*}
\frac{1-e^{-2\lambda-4}}{\lambda+2} = 1,
\end{equation*}
which can be approximated to the machine precision with standard methods (e.g., with MATLAB's \texttt{fzero} we obtain $\lambda = -1.203187869979980$).
The relevant eigenfuction is
\begin{equation*}
\phi(x)=e^{-\int_0^x \mu(s)\dd s-\lambda x} = e^{-(1+\lambda) x}.
\end{equation*}

\begin{figure}
\centering
\includegraphics[trim=28.668pt 0 7.353pt 0]{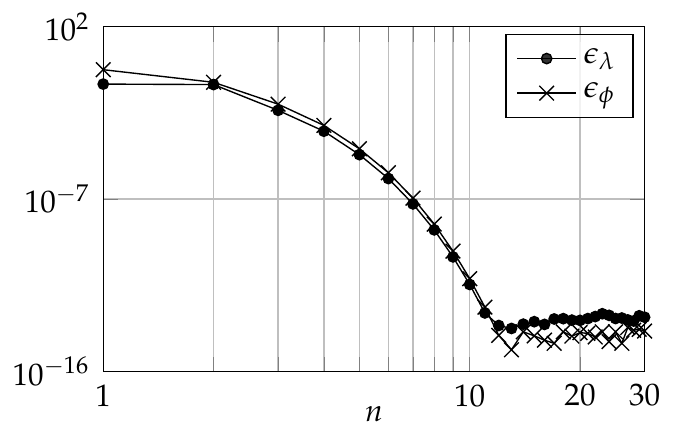}
\caption{Errors $\epsilon_\lambda$ and $\epsilon_\phi$ for the example of appendix \ref{sec:1d}.}
\label{fig:1d}
\end{figure}

We can observe in Figure \ref{fig:1d} that the errors computed by our method decay with infinite order.

\section*{Acknowledgements}
\addcontentsline{toc}{section}{Acknowledgements}

We are grateful to Odo Diekmann for suggesting that we refer to the theory of absolute continuity in the sense of Carathéodory.

The authors are members of INdAM Research group GNCS.
Davide Liessi and Francesca Scarabel are members of UMI Research group ``Mo\-del\-li\-sti\-ca so\-cio-epi\-de\-mio\-lo\-gi\-ca''.
The work of Simone De Reggi and Davide Liessi was supported by the Italian Ministry of University and Research (MUR) through the PRIN 2020 project (No.\ 2020JLWP23) ``Integrated Mathematical Approaches to Socio-Epidemiological Dynamics'' (CUP: E15F21005420006).
Francesca Scarabel is supported by the UKRI through the JUNIPER modelling consortium (grant number MR/V038613/1).

{\sloppy
\printbibliography
\par}

\end{document}